\documentclass[12pt]{amsart}

\usepackage[new]{old-arrows}
\usepackage{amsmath,amssymb,latexsym}
\usepackage{amsfonts,amstext,amsthm,amscd}
\usepackage{graphicx}
\usepackage{enumerate}
\usepackage{amsthm}
\usepackage{xypic}
\theoremstyle{plain}
\usepackage{color}

\usepackage{amssymb}

\makeatletter
\@namedef{subjclassname@2010}{%
  \textup{2010} Mathematics Subject Classification}
\makeatother

\theoremstyle{definition}
\newtheorem{theorem}{Theorem}[section]

\newtheorem{definition}[theorem]{Definition}
\newtheorem{example}[theorem]{Example}
\newtheorem{lemma}[theorem]{Lemma}

\newtheorem{proposition}[theorem]{Proposition}
\newtheorem{remark}[theorem]{Remark}
\newtheorem{question}[theorem]{Question}
\numberwithin{equation}{section}

\numberwithin{equation}{section}

\newcommand{\norm}[1]{\left\Vert#1\right\Vert}
\newcommand{\abs}[1]{\left\vert#1\right\vert}
\newcommand{\Aut}{\mathrm{Aut}}
\newcommand{\B}{\mathfrak B}

\newcommand{\Co}{\mathrm{C}}

\newcommand{\Char}{\mathrm{Char}}

\newcommand{\Ff}{\widehat{f}}
\newcommand{\Homeo}{\mathrm{Homeo}}

\newcommand{\Hom}{\mathrm{Hom}}
\newcommand{\N}{\mathbb{N}}
\newcommand{\Pp}{\mathbb{P}}
\newcommand{\PP}{\mathcal{P}}
\newcommand{\Oo}{\mathcal{O}}
\newcommand{\Ll}{\mathcal{L}}
\newcommand{\Q}{\mathbb{Q}}
\newcommand{\R}{\mathbb{R}}
\newcommand{\Ra}{R_\alpha}
\newcommand{\Rr}{\mathcal{R}}
\newcommand{\Ss}{\mathsf{S}}
\newcommand{\Susp}{\Sigma_f(G)}
\newcommand{\SuspS}{\Sigma_f(\Ss)}

\newcommand{\T}{\mathbb{T}}
\newcommand{\UC}{\mathbb{S}^{1}}
\newcommand{\Z}{\mathbb{Z}}
\newcommand{\Zp}{\Z_p}
\newcommand{\Zz}{\widehat{\Z}}

\newcommand{\To}{\longrightarrow}
\newcommand{\hTo}{\longhookrightarrow}
\newcommand{\mTo}{\longmapsto}
\newcommand{\id}{\mathrm{id}}

\def\bproof{\noindent{\textbf{Proof. }}}
\def\eproof{\noindent{\hfill $\boxdot$}\bigskip}

\begin{document}
\baselineskip=17pt

\title[Some aspects of Rotation Theory on compact abelian groups]{Some aspects of Rotation Theory on compact abelian groups}

\author[M. Cruz-L\'opez]{Manuel Cruz-L\'opez}
\address{ Departamento de Matem\'aticas,
Universidad de Guanajuato, Jalisco s/n, Mineral de Valenciana,
Guanajuato, Gto. 36240 M\'exico.}
\email{manuelcl@ugto.mx}

\author[F. L\'opez-Hern\'andez]{ Francisco J. L\'opez-Hern\'andez}
\address{ Instituto de F\'isica, 
Universidad Aut\'onoma de San Luis Potos\'i, 
Av. Manuel Nava No. 6, Zona Universitaria,
San Luis Potos\'i, SLP. 78290 M\'exico.}
\email{flopez@ifisica.uaslp.mx}

\author[A. Verjovsky]{Alberto Verjovsky}
\address{ Instituto de Matem\'aticas, Unidad Cuernavaca,
Universidad Nacional Aut\'onoma de M\'exico, Apdo. Postal 2 C.P.
2000, Cuernavaca, Mor. M\'exico}
\email{alberto@matcuer.unam.mx}

\begin{abstract}

In this paper we present a generalization of Poincar\'e's Rotation Theory of homeomorphisms of the circle to the case of one-dimensional compact abelian groups which are solenoidal groups, {\it i.e.}, groups which fiber over the circle with fiber a Cantor abelian group. We define rotation elements, \emph{\`a la} Poincar\'e and discuss the dynamical properties of translations on these solenoidal groups. We also study the semiconjugation problem when the rotation element generates a dense subgroup of the solenoidal group. Finally, we comment on the relation between Rotation Theory and entropy for these homeomorphisms, since unlike the case of the circle, for the solenoids considered here there are homeomorphisms (not homotopic to the identity) with positive entropy.\\ 

This is a preliminary version of the published article with minor corrections: Cruz-López, Manuel; López-Hernández, Francisco J.; Verjovsky, Alberto; Some aspects of rotation theory on compact abelian groups. Colloq. Math. 161 (2020), no. 1, 131–155.

\end{abstract}

\thanks{The second author's research was sponsored by CONACYT postdoctoral fellow. The third author was supported by project PAPIIT IN106817, DGAPA, Universidad Nacional 
Aut\'onoma de M\'exico.}
\maketitle

\section[Introduction]{Introduction}
\label{introduction}

Rotation Theory has been an area of intense investigation since its origin in the work of H. Poincar\'e \cite{Poi} in the 1880's.This theory deals with the study of the average behavior of orbits of points in a space under the iteration of functions or flows, but unlike Ergodic Theory, which looks at almost all points, Rotation Theory looks at \emph{all} points. Much of the dynamics is revealed by rotation sets, thus Rotation Theory is part of the Dynamical Systems Theory (see \cite{Mis} for an introductory account of Rotation Theory).

Given $f$ an orientation-preserving homeomorphism of the unit circle $\UC$ Poincar\'e introduced an important invariant, $\rho(f)$, the \emph{rotation number}. This number 
$\rho(f)$ is an invariant of topological conjugacy. Poincar\'e's result in this setting establishes a relationship between the dynamics of $f$ and the translation $R_{\rho(f)}$, whose dynamics depend on the arithmetic nature of $\rho(f)$. More precisely, Poincar\'e's theorem states that for any orientation-preserving homeomorphism $f$ of the unit circle, $\rho(f)$ is rational if and only $f$ has a periodic orbit, and in this case, any nonperiodic orbit is asymptotic to a periodic orbit (see \cite{KH} for a general statement). If the rotation number $\rho(f)$ is irrational then $f$ is semiconjugate to an irrational rotation $R_{\rho(f)}$. The semiconjugacy is actually a conjugacy if the orbits of $f$ are dense. This is a clear illustrationof the relationship between the qualitative behavior of orbits and the rotation set. 
 
The study for homeomorphisms of the $n$\nobreakdash-dimensional torus $\T^n$ is more subtle and many authors have contributed to the theory. For instance, in the 1950s  S. Schwartzman took an interesting approach in the case of an abstract compact metric space using asymptotic cycles \cite{Sch}. 

The basic objective of this work is to describe some aspects of Rotation Theory on general compact abelian topological groups $G$. Our work is motivated and follows the same lines of work done in \cite{CV} for $G=\Ss$, the group called the universal one-dimensional solenoid (see Section \ref{dynamics_solenoid} for a detailed account). Although we discuss the general case, we will continue to focus on one-dimensional solenoidal groups, which are locally compact and Hausdorff abelian groups. 

The rotation set for a homeomorphism $f:G\To G$ isotopic to the identity is a subset of 
$\Hom(\Char(G),\R)$, where $\Char(G)$ is the character group of $G$ consists of continuous homomorphisms from $G$ to $\UC$. From the dynamical perspective, and following the path pioneered by Poincar\'e, the first case to study is that of translations by a specific element of the group. Using a monothetic generator of the group it is possible to describe the topological and ergodic properties of translations. 

In order to get a nontrivial Rotation Theory, one must consider appropriate groups. The following cases are well-established:
\begin{enumerate}
\item If $G$ is the finite group $\Z/n\Z$, the character group of $G$ is isomorphic to itself and therefore $\Hom(\Z/n\Z,\R)=0$, which implies the Rotation Theory is trivial. 

\item When $G$ is the unit circle $\UC$, $\Char(\UC)\cong \Z$ and therefore 
$\Hom(\Z,\R)\cong \R$. The rotation set consists of a single point: the usual translation number of Poincar\'e \cite{Poi}, \cite{Ghys}, \cite{Her}.

\item For $G$ the 2-dimensional torus $\T^2$, $\Char(\T^2)\cong \Z\times \Z$ and therefore $\Hom(\Z\times\Z)\cong \R\times \R$. The rotation set is the usual translation set \cite{Fra}, \cite{Jag}, \cite{MZ}.
\end{enumerate}

Any compact abelian torsion-free group $G$ admits a decomposition of the form
\[ \bigoplus_{r_0} \Ss \times \prod_{p\in \Pp} \left( \bigoplus_{r_p} \Zp \right) \]
where $r_0$ and $r_p$ are cardinal numbers, $\Ss$  is the universal one-dimensional solenoidal group see  Example \ref{universal_solenoid} for the full definition), $\Zp$ is the Cantor group of $p$-\emph{adic} integers and $\Pp$ is the set of all prime numbers (see \emph{e.g.}, \cite{HR,Pon}). 

For these more general groups what is currently known about Rotation Theory is as follows:

\begin{enumerate}
\item If $G$ is a totally disconnected group, the isotopy component of the identity in the group of homeomorphisms is trivial.
\item If $G$ is a one-dimensional solenoidal group (\emph{i.e.}, the Pontryagin dual of a dense subgroup of $\Q$ with the discrete topology) initial steps were developed in \cite{CV}. These results are described in Section \ref{dynamics_solenoid}.
\end{enumerate}

For a general compact abelian group $G$ there is no such structure theorem, even though some possible factors appearing in the decomposition mentioned above are groups already discussed. From the viewpoint of Rotation Theory, it is not enough to study each factor separately, as illustrated by the Rotation Theory of the torus $\T^2$.

Finally, another important aspect analyzed in this work is the relationship between Rotation Theory and entropy (see a historical account of entropy in dynamics in \cite{Kat}). This topic has been the subject of intense study during the last 20 years for instance in \cite{Cro,KH,KW1,KW2,Kwa2,LlM}. 


The article is organized as follows. Section \ref{compact_abelian-groups} introduces general aspects of compact abelian groups including duality theory. Section 
\ref{homeo-susp_G}  describes the suspension of a homeomorphism of a compact abelian group which is isotopic to the identity. Section \ref{rotation_theory} then describes rotation sets for homeomorphisms of compact abelian groups which are isotopic to the identity (this is done in two slightly different ways derived from asymptotic cycles). Dynamical properties of translations are presented in Section \ref{dynamics_translations} and Section \ref{basic_examples} is a summary of Poincar\'e's Rotation Theory on three fundamental examples: the circle, the torus and the universal one-dimensional solenoid. Section \ref{entropy_rotation} summarizes the relationship between entropy and rotation sets as described in several recent works. Finally, in Section \ref{final_remarks}, we comment on the group of homeomorphisms of a general compact abelian group.

\section[Compact abelian groups]{Compact abelian groups}
\label{compact_abelian-groups}

This section is a brief summary of relevant portions of the theory of compact abelian groups. It provides basic definitions and introduces the character group and Pontryagin duality.

Recall we are focusing on topological abelian groups which are locally compact and Hausdorff. 

\subsection[Abelian groups]{Abelian groups}
\label{abelian_groups}

Let $G$ be a set with a binary operation $\cdot:G\times G\To G$ and let $\tau$ be a family of subsets of $G$. The set $G$ is called a \emph{topological group} if
\begin{enumerate}
\item $(G,\cdot)$ is a group,

\item $(G,\tau)$ is a topological space,

\item The maps $m:G\times G\To G$ and $i:G\To G$ given by
\[ (x,y)\mTo x\cdot y, \qquad x\mTo x^{-1} \]
respectively are continuous.
\end{enumerate}
The operation of the group will be denoted by $+$ and the identity element by 0. 

Here are some general examples:

\begin{example}
\label{examples1}
\begin{enumerate}
\item Any abelian group equipped with the discrete topology is a topological abelian group.
\item If $\{G_\alpha\}_{\alpha\in \Lambda}$ is a family of topological abelian groups, the product space $\displaystyle{\prod_{\alpha\in \Lambda} G_\alpha}$ with the natural group operation and the product topology is a topological abelian group.
\item A projective limit $\displaystyle{\varprojlim_{_{\alpha,\beta\in \Lambda}}  
\left\{ \varphi_{\alpha\beta} : G_\beta\To G_\alpha \right\}}$ obtained from a directed set 
$\Lambda$, where $\left\{\varphi_{\alpha\beta} \right\}_{_{\alpha,\beta\in\Lambda}}$ are epimorphisms of compact abelian groups, is a compact abelian group. 
\end{enumerate}
\end{example}

An important special instance of Example \ref{examples1}(3) above is the following. 

\begin{example} [Profinite completion of $\Z$]
\label{profinite_completionZ}
The \emph{profinite completion} of $\Z$ is
\begin{eqnarray*}
\Zz & = & \varprojlim_{n\in \N} \Z/n\Z \\
 & =&  
\left\{ x=(a_n)_{n\in \N} \in \prod_{n\in \N} \Z/n\Z \,: \,
a_m \equiv a_n \,\text{mod}\, m, \,\; \forall\; m\mid n \right\}.
\end{eqnarray*}

With the group structure and topology induced from the product space, $\Zz$ is a compact abelian and Hausdorff topological group. With respect to this profinite topology, $\Zz$ is totally disconnected and perfect, and hence homeomorphic to the Cantor set \cite{Wil}. Moreover, $\Zz$ admits a dense inclusion of $\Z$.
\end{example}

Let $G$ be a locally compact abelian group. If $g\in G$, the following maps are well-defined:
\[ \text{(Left translation)}\qquad L_g:G\To G, \quad x\mTo x-g \]
\[ \text{(Right translation)}\qquad R_g:G\To G, \quad x\mTo x+g \]
Recall that we had previously defined 
\[ \text{(Inverse Map)} \qquad i:G\To G\qquad x\mTo x^{-1}\]

\begin{proposition} 
The maps $L_g$, $R_g$ and $i$ are homeomorphisms.
\end{proposition}

If $H$ is a subgroup of $G$, the set of cosets of $H$ in $G$ is called the \emph{quotient space} of $G$ by $H$, denoted by $G/H$. The canonical projection 
\[ \pi:G\To G/H \text{ is given by } x\mTo [x]:=xH. \] 
In $G/H$ we define the quotient topology as the finest topology which makes the projection map continuous. That is, $U\subset G/H$ is open if and only if $\pi^{-1}(U)\subset G$ is open.

\begin{remark} 
If $H$ is a subgroup of $G$ then $G/H$ is Hausdorff if and only if $H$ is closed.
\end{remark}

\begin{example} [Circle]
\label{circle} 
The quotient group $\R/\Z$ is called the one-dimensional torus and it is denoted by $\T$. This abelian group is isomorphic to the group of the circle $\UC$ under the exponential map. 
\end{example}

\begin{example} [Tori]
\label{n-torus}
The $n$-dimensional torus 
\[ \T^n =\underbrace{\UC \times \cdots \times \UC}_{\; \text{n times}} \] 
and the infinite torus 
\[ \T^{\infty} = \prod_{n\geq 0} \UC \] 
are compact abelian topological groups.
\end{example}

\begin{example} [Universal Solenoid]
\label{universal_solenoid} 
The subset of $\T^{\infty}$
\[ \Ss := \{ (x_n)_{n\geq 0} \in \T^{\infty} : x_m^{m-n}=x_n,\; \rm{when} \; n\mid m\} \]
is a compact abelian topological group called the \emph{universal one-dimensional solenoid}. The projection of $\Ss$ onto the first factor determines a $\Zz$-bundle structure 
$\Zz\hTo \Ss\To \UC$. Furthermore, $\Ss$ admits a dense one-parameter subgroup 
$\sigma:\R\hTo \Ss$ and it is locally homeomorphic to the product of the Cantor set and an interval. Because $\Ss$ contains a dense one-parameter subgroup, $\Ss$ is called a \emph{solenoidal group}. Denote by $\Ll_0$ the image of $\R$ under $\sigma$. One has that $\Ll_0$ is the connected component of the identity.
\end{example} 

\subsection[Characters and duality]{Characters and duality}
\label{char-duality_G}

The set 
\[ \Char(G) := \Hom_{\text{cont}}(G,\UC) \] 
consisting of all continuous homomorphisms from $G$ into the multiplicative group 
$\UC$, is a topological abelian group with respect to the operation
$$ 
(\chi_1 \cdot \chi_2)(x) = \chi_1(x) \chi_2(x),\quad (x\in G,\, \chi_1,\chi_2\in \Char(G)) 
$$
and the uniform topology. This group is called the \emph{character group}, or the \emph{Pontryagin dual} of $G$. 

The character group has the following properties.

\begin{remark} 
\begin{enumerate}
\item If $G$ is a compact abelian group then $\Char(G)$ is a discrete abelian group. 
\item If $G$ is a discrete abelian group then $\Char(G)$ is a compact abelian group.
\item (Pontryagin duality) $\Char(\Char(G)) \cong G$.
\end{enumerate}
\end{remark}

The following are important examples:
\begin{example}
\begin{enumerate}
\item $\Char(\UC) \cong \Z$.
\item $\Char(\T^n) \cong \Z^n$.
\item $\Char(\Zz) \cong \Q/\Z$, where $\Q/\Z$ is the group of roots of unity.
\item $\Char(\Ss) \cong  \Q$.
\end{enumerate}
\end{example}

\begin{example}[Solenoids]
1-dimensional compact abelian solenoidal groups are Pontryagin duals of subgroups of the additive rational numbers $\Q$ with the discrete topology. They are inverse limits of directed sets of epimorphisms of the circle onto itself (see Example \ref{universal_solenoid}).
\end{example}

\begin{remark}
\cite{Ste} Suppose $G$ is also compact and connected. If $\check{H}^1(G,\Z)$ denotes the first \v{C}ech cohomology group of $G$ with coefficients in $\Z$ then,
\[ \check{H}^1(G,\Z)\cong \Char(G).\]
\end{remark}

\section[Homeomorphisms and suspensions]{Homeomorphisms and suspensions}
\label{homeo-susp_G}

We describe the suspension of a homeomorphism of a general compact abelian group and focus on homeomorphisms isotopic to the identity. Some comments on the character group of a suspension space are made and measures defined on groups and suspensions are also considered. 

The suspension of a homeomorphism is a special case of the much general case of the suspension of a representation. This construction fits in the context of Foliation Theory. A very complete treatment on the subject of Foliation Theory is \cite{CC}. The character group of the suspension is treated via the Bruschlinsky-Eilenberg's theory as in \cite{Sch}. 

\subsection[The suspension of a homeomorphism]{The suspension of a homeomorphism}
\label{suspension_homeomorphism}

Denote by $\Homeo(G)$ the group of all homeomorphisms $f:G\To G$. Suppose that 
$f\in\Homeo(G)$ and consider on $G\times [0,1]$ the equivalence relation
\[ (z,1)\sim (f(z),0)\qquad (z\in G). \]

The \emph{suspension} of $f$ is the compact space
\[ \Susp := G\times [0,1] /(z,1)\sim (f(z),0). \]
In $\Susp$ there is a well-defined flow 
$\phi:\R\times \Susp\To \Susp$, called the \emph{suspension flow} of $f$, given by
\[ \phi(t,[(z,s)]) :=\phi_t([(z,s)]) = [(f^m(z),t+s+m)], \] 
if $m\leq t+s < m+1$. Here $f^m = \underbrace{f\circ \cdots \circ f}_{m\; \text{times}}$
denotes the $m$th iteration of $f$.

The canonical projection $\pi:G\times [0,1]\To \Susp$ sends $G\times \{0\}$ homeomorphically onto its image $\pi(G\times \{0\})\equiv G$ and every orbit of the suspension flow intersects $G$. The orbit of any $(z,0)\in \Susp$ must coincide with the orbit $\phi_t(z,0)$ at time $0\leq t\leq T$ for $T$ an integer.

The subgroup of $\Homeo(G)$ containing all the homeomorphisms which are isotopic to the identity will be denoted by $\Homeo_+(G)$. When $f\in\Homeo_+(G)$, the suspension 
$\Susp$ is homeomorphic to $G\times \UC$ hence if $G$ is a compact abelian group, 
$\Susp$ is also a compact abelian group.  

\begin{remark} 
The compact abelian topological group structure of $G\times \UC$ will play a fundamental role in the writing. In what follows when $f\in\Homeo_+(G)$ we will abuse the notation by considering $\Susp$ as a compact abelian topological group with the structure of 
$G\times \UC$. 
\end{remark}

\subsection[The Universal Solenoid as the suspension of a translation in $\Zz$ ]{The Universal Solenoid as the suspension of a translation in $\Zz$ }
\label{suspension_representation}

We illustrate the general construction of the suspension of a representation by describing the universal solenoid as the foliated space obtained from the suspension process. Recall that the fundamental group of the circle $\pi_1(\UC)$ is isomorphic to the group of integers $\Z$. 
Also, $\Zz$ admits a dense inclusion $i:\Z\to\Zz\,$ of $\Z$. Let $i(n)=\mathbf{n}$.
Denote by $\Homeo(\Zz)$ the group consisting of all homeomorphisms of the Cantor group 
$\Zz$. Take a representation $\pi_1(\UC)\To \Homeo(\Zz)$ acting by translations and consider the action of $\Z$ on $\R\times \Zz$:
\[ \Z \times (\R\times \Zz) \To \R\times \Zz, \quad (n,(t,x))\mTo (t+n,x+\mathbf{n}), \, n\in\Z,t\in\R,x\in\Zz\]

Here, $\R$ is thought of as the universal covering space of the circle $\UC$. The action of 
$\Z$ in the first component is by deck transformations and in the second coordinate is by translations inside the profinite completion $\Zz$. The universal solenoid $\Ss$ is the foliated space obtained as the quotient $\R\times_{\Z} \Zz$.

\begin{example}
If $f\in\Homeo_+(\UC)$ then the suspension space $\Sigma_f(\UC)$ is homeomorphic to the torus $\T^2$. The suspension flow is 
\[ \phi(t,(z,s)) = (f^m(z),t+s+m), \qquad t\in\R, \; (z,s)\in\T^2, \] 
if $m\leq t+s < m+1$ (see Figure 1).
\begin{center}
\includegraphics[scale=0.25]{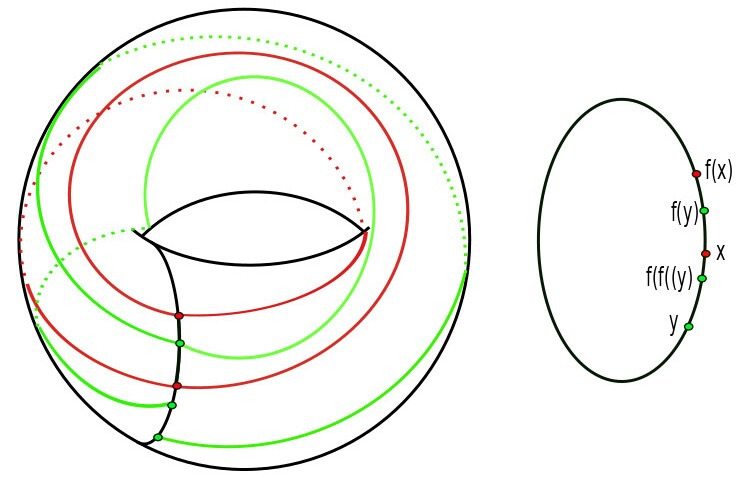}
\centerline{Figure 1. Suspension flow with a periodic orbit of period 2}
\end{center}
\end{example}

 The universal solenoid $\Ss$ is the suspension space of a homeomorphism of the Cantor group $\Zz$. In fact, it is the suspension of translation by $\mathbf1$ in $\Zz$ (adding machine).

\begin{example} [The solenoid as a suspension]
If $R_1:\Zz\To \Zz$ is the homeomorphism of translation by 
$\mathbf1=i(1)\,\, \text{(}i:\Z\to\Zz \, \text{is the canonical inclusion)}$ them the suspension space
\[ \Sigma_{R_1} (\Zz)= \Zz\times [0,1] /(s,1)\sim (R_1(s),0) \]
is homeomorphic to the solenoid $\Ss$. This can be explained as follows: Identify any element $s$ in the canonical fiber $\Zz\hTo \Ss$ with $[(s,1)]\in \Sigma_{R_1}(\Zz)$, and any element $s+x\in\Ss,\; (x\in (0,1)\subset \Ll_0)$ with the element $[(s,x)]\in \Sigma_{R_1}(\Zz)$. The one-parameter subgroup $\sigma:\R\hTo \Ss$ defines a flow 
$\R\times \Ss\To \Ss$: $(t,z)\mTo \sigma(t) + z$. Clearly this flow coincides with the suspension flow $\phi_t$.

Note that $R_1$ is not isotopic to the identity, so the topological structure of 
$\Sigma_{R_1}(\Zz)$ is not the topological product structure $\Zz\times \UC$, in fact 
$\Sigma_{R_1}(\Zz)$ is connected.
\end{example}

\subsection[Characters and suspensions]{Characters and suspensions}
\label{characters_suspensions}

Denote by $\Co(\Susp,\UC)$ the topological abelian group of continuous functions from  
$\Susp$ to $\UC$ under pointwise multiplication and with the compact-open topology. 

The subset $R(\Susp,\UC)\subset \Co(\Susp,\UC)$ consisting of continuous functions 
$h:\Susp\To \UC$ that can be written as $h(z,s) = \exp(2\pi i\psi(z,s))$ with 
$\psi:\Susp\To \R$ a continuous function, is a closed subgroup. Hence the quotient group 
$\Co(\Susp,\UC)/R(\Susp,\UC)$ is a topological group. By Bruschlinsky-Eilenberg's theory,
\[ \check{H}^1(\Susp,\Z)\cong \Co(\Susp,\UC)/R(\Susp,\UC). \]

Since
\[ \check{H}^1(\Susp,\Z)\cong \Char(\Susp), \] 
it follows that
\[ \Char(\Susp)\cong \Co(\Susp,\UC)/R(\Susp,\UC). \]

On the other hand, since $\Susp$ has the same topological product group structure as 
$G\times \UC$, this implies that its character group is given by
\[ \Char(\Susp)\cong \Char(G)\times \Char(\UC)\cong \Char(G)\times \Z. \]

Given any element $(k_g,n)\in \Char(G)\times \Z$, the corresponding character 
$\chi_{k_g,n}\in \Char(\Susp)$ can be written as
\begin{align*}
\chi_{k_g,n}(z,s)	&= \exp(2\pi ik_g z)\cdot \exp(2\pi ins)\\
							&= \exp(2\pi i (k_g z + ns))
\end{align*}
for any $(z,s)\in \Susp$, where $k_g z$ means $z$ added with itself $k_g$ times.

\subsection[Measures]{Measures}
\label{measures}

Invariant Borel probability measures play an important role in our definitions of rotation elements and rotation sets as well as in ergodic theory of dynamical systems. Given any $f$-invariant Borel probability measure $\mu$ on $G$ and $\lambda$ the usual Lebesgue measure on $[0,1]$, the product measure $\mu\times \lambda$ allows to define a $\phi_t$-invariant Borel probability measure on $\Susp$. Reciprocally, given any $\phi_t$-invariant Borel probability measure $\nu$ on $\Susp$, it can be defined by disintegration with respect to the fibers, an $f$-invariant Borel probability measure $\mu$ on $G$. Denote by $\PP(G)$ the weak$^*$ compact convex space of $f$-invariant Borel probability measures defined on $G$.

\section[Rotation Theory on compact abelian groups]{Rotation Theory on compact abelian groups}
\label{rotation_theory}

General aspects of Rotation Theory for compact abelian groups are treated in this section using the notion of 1-cocycles and asymptotic cycles. 

\subsection[1-cocycles]{1-cocycles}
\label{1-cocycles}

A $1$\nobreakdash-\emph{cocycle} associated to the suspension flow $\phi_t$ is a continuous function
\[ C:\R\times \Susp\To \R \]
which satisfies the relation
\[ C(t+u,(z,s)) = C(u,\phi_t(z,s)) + C(t,(z,s)), \] 
for every $t,u\in \R$ and $(z,s)\in \Susp$. The set which consists of all $1$\nobreakdash-cocycles associated to $\phi_t$ is an abelian group denoted by $\Co^1(\phi)$. A 1-\emph{coboundary} is the 1-cocyle determined by a continuous function 
$\psi:\Susp\To \R$ such that
\[ C(t,(z,s)) := \psi(z,s) - \psi(\phi_t(z,s)). \] 

The set of $1$\nobreakdash-coboundaries $\Gamma^1(\phi)$ is a subgroup of 
$\Co^1(\phi)$ and the quotient group
\[ H^1(\phi):=\Co^1(\phi)/\Gamma^1(\phi) \]
is called the $1$\nobreakdash-\emph{cohomology group} associated to $\phi_t$. The proof of the next Proposition (for an arbitrary compact metric space) can be seen in \cite{Ath}.

\begin{proposition}
\label{associated_cocycle} 
For every continuous function $h:\Susp\To \UC$ there exists a unique $1$\nobreakdash-cocycle 
$C_h:\R\times \Susp\To \R$ associated to $h$ such that
\[ h(\phi_t(z,s)) = \exp(2\pi iC_h(t,(z,s)))\cdot h(z,s), \] 
for every $(z,s)\in \Susp$ and $t\in \R$.
\end{proposition}

This Proposition implies that there is a well-defined homomorphism
\[ \Char(\Susp)\cong \check{H}^1(\Susp,\Z)\To H^1(\phi) \] 
sending any character $\chi_{k_g,n}\in \Char(\Susp)$ to the cohomology class 
$[C_{\chi_{k_g,n}}]$, where $C_{\chi_{k_g,n}}$ is the unique $1$\nobreakdash-cocycle associated to $\chi_{k_g,n}$.

Applying the above Proposition to any nontrivial character $\chi_{k_g,n}\in \Char(\Susp)$ the following relation is obtained:
\[ \chi_{k_g,n}(\phi_t(z,s)) = \exp(2\pi iC_{\chi_{k_g,n}}(t,(z,s)))\cdot \chi_{k_g,n}(z,s). \] 

Using the explicit expressions for the characters on both sides of the above equation, the next equalities hold
\begin{align*}
\chi_{k_g,n}(\phi_t(z,s))	&= \chi_{k_g,n}(f^m(z),t+s-m)\\
										&= \exp(2\pi i(k_gf^m(z) + n(t+s-m)))\\
										&= \exp(2\pi i(k_gf^m(z) + nt + ns))
\end{align*} and

\[ \chi_{k_g,n}(z,s) = \exp(2\pi i(k_gz+ns)). \]

Comparing these two expressions one gets
\begin{equation}
\label{cocycle1} 
C_{\chi_{k_g,n}}(t,(z,s)) = k_g(f^m(z)-z) + nt.
\end{equation}

In the following sections we will be using Birkhoff's Ergodic Theorem to study rotation sets, therefore the computation at time $t=1$ is required. If $t=1$ then $m=1$ and the 
$1$\nobreakdash-cocycle at time $t=1$ is
\begin{equation}
\label{time1cocycle} 
C_{\chi_{k_g,n}}(1,(z,s)) = k_g(f(z)-z) + n.
\end{equation}

\subsection[Rotation sets for homeomorphisms isotopic to the identity]{Rotation sets for homeomorphisms isotopic to the identity}
\label{rotation_element-isotopic}

\subsubsection{Rotation sets through cocycles}

Suppose that $f:G\To G$ is a homeomorphism isotopic to the identity. If $\nu$ is any 
$\phi_t$-invariant Borel probability measure on $\Susp$, by Birkhoff's ergodic theorem there is a well-defined homomorphism $H^1(\phi)\To \R$ given by
\[ [C_\chi]\mTo \int_{\Susp} C_\chi(1,(z,s)) d\nu. \]

Composing the two homomorphisms
\[ \Char(\Susp)\To H^1(\phi)\To \R \] 
it is obtained a well-defined homomorphism $H_{f,\nu}:\Char(\Susp)\To \R$ given by
\[ H_{f,\nu}(\chi_{k_g,n}) := \int_{\Susp} C_{\chi_{k_g,n}}(1,(z,s)) d\nu. \]

Denote by $\mu$ the $f$-invariant Borel probability measure on $G$ obtained by disintegration of $\nu$ with respect to the fibers. Evaluating the above integral using equation (\ref{time1cocycle}) yields
\begin{align*}
H_{f,\nu}(\chi_{k_g,n})		&= \int_{\Susp} (k_g(f(z)-z)  + n)d\nu\\
										&= k_g\int_G (f(z)-z)  d\mu + n.
\end{align*}

Hence $H_{f,\nu}$ determines an element in $\Hom(\Char(\Susp),\R)$ for each measure 
$\nu$ in $\Susp$, and therefore, for each measure $\mu\in \PP(G)$. Hence one gets a well-defined function
\[ H_f:\PP(G)\To \Hom(\Char(\Susp),\R) \] 
given by $\mu\mTo H_{f,\mu}$, where $H_{f,\mu}$ is:
\[ H_{f,\mu}(\chi_{k_g,n}) = k_g\int_G (f(z)-z) d\mu + n. \]

By composing $H_f$ with the continuous homomorphism
\[ \Hom(\Char(\Susp),\R)\To \Char(\Char(\Susp)) \] 
given by
\[ H_{f,\mu}\mTo \pi\circ H_{f,\mu}, \] 
where $\pi:\R\To \UC$ is the universal covering projection, we obtain a well-defined continuous function $\rho:\PP(G)\To \Char(\Char(\Susp))$ given by
\[ \mu\mTo \rho_\mu := \pi\circ H_{f,\mu}. \]

The above discussion is summarized in the following theorem:

\begin{theorem}
\label{homrot}
For each $\mu\in \PP(G)$, there exists a well-defined continuous homomorphism 
\[ \rho_\mu:\Char(\Susp)\To \UC \] 
given by
\begin{align*}
\rho_\mu(\chi_{k_g,n})	&:= \exp(2\pi iH_{f,\mu}(\chi_{k_g,n}))\\
										&= \exp \left(2\pi ik_g \int_G (f(z)-z) d\mu\right).\\
\end{align*}
\end{theorem}

By Pontryagin's duality theorem,
\[ \Char(\Char(\Susp))\cong \Susp \]
and therefore $\rho_\mu\in \Susp$. Since $\Susp\simeq G\times \UC$ and 
$\rho_\mu(\chi_{k_g,n})=\rho_\mu(\chi_{k_g,0})$, it follows that $\rho_\mu$ does not depend on the second component so we can identify $\rho_\mu$ with 
$(\rho_\mu,1)\in G\times \UC$.

\begin{definition}
The element $\rho_\mu(f) := \rho_\mu \in G$ defined above is the \emph{rotation element} associated to $f$ with respect to the measure $\mu$.
\end{definition}

\begin{remark} 
By definition, $\rho_\mu(f)$ can be identified with the element $\int_G (f(z)-z) d\mu$ in 
$G$ determined by the character of $G$ given by 
\[ k_g\mTo \exp \left(2\pi ik_g\int_G (f(z)-z) d\mu \right). \] 

That is, $\rho_\mu(f)$ is \emph{G-valued}.
\end{remark}

If $\rho:\PP(G)\To G$ is the map given by $\mu\mTo \rho_\mu(f)$, then 
$\rho$ is continuous from $\PP(G)$ to $G$. Since $\PP(G)$ is compact and convex, and $f$ is isotopic to the identity, the image $\rho(\PP(G))$ is a compact subset of $G$. 

\begin{definition} [Rotation set]
The set $\rho(\PP(G))$ is the \emph{rotation set} of $f$.
\end{definition}

\subsubsection[Rotation sets through functionals]{Rotation sets through functionals}
\label{rotation_element-functionals}

Recall that the group \\ $\Hom(\check{H}^{1}(\Susp,\Z),\R)$ consists of all homomorphisms from $\check{H}^{1}(\Susp,\Z)$ to $\R$. Note:
\begin{align*}
\Hom(\check{H}^{1}(\Susp,\Z),\R)	&\cong \Hom(\Char(G) \times\Z,\R) \\
&\cong \Hom(\Char(G),\R)\times \Hom(\Z,\R) \\
&\cong \Hom(\Char(G),\R)\times \R.
\end{align*}

Given $g\in \Susp$ and $T>0$ there is a continuous linear functional $\Lambda_{g,T}$ in 
$\Hom(\Co(\Susp,\UC),\R)$ defined by
\[ \Lambda_{g,T}(k) = \int_{0}^{T} \frac{d}{dt} \arg[k](\phi_t(g)) dt, \]
where $k:\Susp\To \UC$ is a continuous function with argument $\arg[k]$ such that the argument is a continuous function along the orbits of the suspension flow. These linear functionals can be described by the following result of Schwartzman (see \cite{Pol}).

\begin{proposition} 
For each $g\in \Susp$ the following conditions hold:
\begin{enumerate}
\item The family of linear functionals 
$\{\Lambda_{g,T}\}_{T\in \R_+}\subset \Hom(\Co(\Susp,\UC),\R)$ is equicontinuous in the weak* topology.
\item The limit points $\Rr^g\subset \Hom(\Co(\Susp,\UC),\R)$ of this family are constant on equivalence classes of maps in $\Co(\Susp,\UC)$.
\end{enumerate}
\end{proposition}

This Proposition tells us that each $r^g\in \Rr^g$ gives a well-defined element in 
$\Hom(\check{H}^{1}(\Susp,\Z),\R)$ and therefore, using the comments above, an element in $\Hom(\Char(G),\R)\times \R$. As a consequence the linear functionals can be written in terms of their coordinates 
\[ r^g = (r^g_1,r^g_2)\in \Rr^g\subset \Hom(\Char(G),\R)\times \R. \] 

Due to the explicit construction of the functional $\Lambda_{g,T}$ and $\Rr^g$, it is easily seen that $r^g_2\equiv 1\in \R$. Therefore $r^g$ is independent of the second coordinate in $\Susp$. Furthermore, if $g=(z,s)$ then $r^g_1$ depends only on the first coordinate 
$z\in G$. Therefore there is only need to consider the components $r^g_1\in \R$ of the linear functionals $r^g\in \Rr^g$. Let us denote $r^g_1$ by $r^z$.

\begin{definition}
Suppose $f\in \Homeo_+(G)$ and $z\in G$. The \emph{punctual rotation set} of $f$ relative to $z$ is 
\[ \rho_z(f) = \{ r^z : r^g\in \Rr^g \}\subset \Hom(\Char(G),R). \]

The \emph{rotation set} of $f$ is defined as
\[ \rho(f) = \bigcup_{z\in G} \rho_z(f). \]
\end{definition}

\subsection[The rotation element \emph{\`a la} de Rham]{The rotation element \emph{\`a la} de Rham}
\label{rotation_deRham}

If $d\lambda$ denotes the usual Lebesgue measure on $\UC$, then, given any character 
$\chi_{k_g,n}\in \Char(\Susp)$ there is a well-defined closed differential one form on 
$\Susp$ given by
\[ \omega_{\chi_{k_g,n}} = \chi_{k_g,n}^* d\lambda. \]

Let $X$ be the vector field tangent to the flow $\phi_t$ and let $\nu$ be any $\phi_t$-invariant Borel probability measure on $\Susp$. Define
\[ H_{f,\nu}:\Char(\Susp)\To \R \]
by
\[ H_{f,\nu}(\chi_{k_g,n}) = \int_{\Susp} \omega_{\chi_{k_g,n}}(X) d\nu \]
and observe that this definition only depends on the cohomology class of 
$\omega_{\chi_{k_g,n}}$ and the measure class of $\nu$. Hence, there is a well-defined continuous homomorphism 
$\rho(f):\Char(\Susp)\To \UC$ given by
\[ \rho(f)(\chi_{k_g,n}) = \exp(2\pi iH_{f,\nu}(\chi_{k_g,n})). \]

Thus, as before,
\[ \rho(f)\in \Char(\Char(\Susp))\cong \Susp. \]

\begin{proposition} 
$\rho(f)$ is the rotation element associated to $f$ corresponding to $\nu$.
\end{proposition}

\subsection[Rotation sets for general homeomorphisms]{Rotation sets for general homeomorphisms}
\label{rotation_element-homeomorphism}

For a general homeomorphism $f\in \Homeo(G)$, the construction made in Section \ref{rotation_element-isotopic} to define the rotation element can be rewritten, \emph{mutatis mutandis}, to get a rotation element in the following way:
\begin{enumerate}[a.]
\item The suspension of $f$ is the compact space
\[ \Susp = G\times [0,1] /(z,1)\sim (f(z),0). \]

\item The suspension flow of $f$ is well-defined and given by
\[ \phi(t,[z,s]) = [f^m(z),t+s-m], \] 
if $m\leq t+s < m+1$.

\item There are isomorphisms
\[ \Char(\Susp)\cong \check{H}^1(\Susp,\Z)\cong \Co(\Susp,\UC)/R(\Susp,\UC). \]

\item The image of the well-defined homomorphism
\[ H_{f,\nu} : \Char(\Susp)\To \R, \qquad \chi\mTo \int_{\Susp} C_{\chi}(1,\cdot) d\nu, \] 
is the set (see \cite{Ath})
\[ \left\{ \int_G \psi d\mu : \psi\in \mathrm{Log}(G,f) \right\}, \] 
where
\[ \mathrm{Log}(G,f) = \{\psi:G\To \R: h\circ f = \exp(2\pi i\psi)h,\; h\in \Co(G,\UC)\}. \]

\item Define
\[ \rho(f)_\mu([h]) := \exp \left( 2\pi i\int_G \psi d\mu \right). \]
\end{enumerate}

\section[Dynamics of translations on compact abelian groups]{Dynamics of translations on compact abelian groups}
\label{dynamics_translations}

This section reviews the topological dynamics of translations on compact abelian groups. This classical theory can be consulted e.g. \cite{CFS,KH}.

Let $G$ be a compact abelian group and let $\Ra$ be the translation by $\alpha\in{G}$.
Recall that the \emph{orbit} of $z$ under $\Ra$ is 
\[ \Oo_{\Ra}(z) = \{\Ra^n(z): n\in \Z\}. \]

The basic problem in iteration theory is:\\

\noindent \textsc{Problem:} To describe the asymptotic behavior of orbits of points in $G$ under iteration of $\Ra$.

The following general definition will be useful in our analysis.

\begin{definition} 
An element $\alpha\in G$ is a \emph{monothetic generator} of $G$ if the set of translates of $\alpha$, $T(\alpha)=\{n\alpha:n\in \Z\}$ is dense in $G$. If there exists a monothetic generator of $G$, $G$ is called a \emph{monothetic} group.
\end{definition}

\begin{remark} 
If $G$ is a monothetic group with monothetic generator $\alpha$, then, since $T(\alpha)$ is a (cyclic) commutative group, its closure $G$ is also abelian. Therefore any monothetic group is automatically abelian.
\end{remark}

Observe that the orbit of 0 under the translation $R_\alpha$ coincides with the set 
$T(\alpha)$. The following examples illustrates this notion.

\begin{example}
\begin{enumerate}
\item (Cantor group $\Zz$) As a consequence of Example \ref{profinite_completionZ}, the orbit of $\mathbf{0}\in \Zz$ under the translation $R_{\mathbf{1}}$ is dense in $\Zz$. Since this orbit coincides with $T(\mathbf{1})$, this implies that $\Zz$ is a monothetic group with monothetic generator $\mathbf{1}$.
\item (The circle group) Si $\alpha\in \UC$ is irrational then the orbit of $0\in \UC$ under the translation $R_\alpha$ is dense. Hence $\UC$ is a monothetic group with monothetic generator $\alpha$.
\item (Tori) Si $\alpha\in \T^n$ is totally irrational then the orbit of $0\in \T^n$ under the translation $R_\alpha$ is dense. Hence $\T^n$ is a monothetic group with monothetic generator $\alpha$.
\item (Universal solenoid) $\Ss$ is a monothetic group as showed in \cite{FHK}, Section 2. 
\end{enumerate}
\end{example}

\subsection[Minimality of translations]{Minimality of translations}
\label{minimality_translations}

A homeomorphism $f:X\To X$ of the topological space $X$ is said to be \emph{minimal} if every orbit is dense.

\begin{proposition} 
If the translation $\Ra:G\To G$ has a dense orbit, then $\Ra$ is minimal.
\end{proposition}

\bproof 
Let $z,w\in G$ be any two points with orbits $\Oo_{\Ra}(z)$ and $\Oo_{\Ra}(w)$, respectively. Hence
\begin{align*}
\Ra^n(w)	&= w + n\alpha\\
				&= z + n\alpha + (w-z)\\
				&= \Ra^n(z) + (w-z).
\end{align*} 

Therefore
$$
\overline{\Oo_{\Ra}(z)}=G \qquad \text{if and only if}\qquad \overline{\Oo_{\Ra}(w)}=G.
$$
\eproof

The following result is clear:

\begin{theorem} 
If $\alpha\in G$, the following conditions are equivalent:
\begin{enumerate}[a.]
\item The translation $\Ra:G\To G$ is minimal.
\item $G$ is a monothetic group with generator $\alpha$.
\end{enumerate}
\end{theorem}

\subsection[Ergodicity of translations]{Ergodicity of translations}
\label{ergodicity_translations-groups}

This section analyzes ergodic properties of translations. First, some general definitions on measures and integration are briefly introduced.

\subsubsection[Haar measure]{Haar measure}
\label{haar_measure} 

For completeness, let's recall the definition and some of the properties of the Haar measure of a locally compact abelian group.

\begin{definition}
A \emph{right Haar measure} in $G$ is a Borel measure $\mu$ in $G$ which satisfies the following properties:
\begin{enumerate}[$\bullet$]
\item $\mu$ is finite on compact subsets, in particular, $\mu(G)<\infty$.

\item $\mu$ is regular in all Borel subsets.

\item $\mu(U)>0$ for each open subset $U$ of $G$.

\item $\mu$ is invariant under right translations:
\[ \mu(\Ra^{-1}(E)) = \mu(E)\qquad \text{for all}\; E\in \B(G). \]
\end{enumerate}
\end{definition}

Since $G$ is compact we can normalize the Haar measure on $G$ such that $\mu(G)=1$. The Haar measure always exists on any compact abelian topological group.

\begin{theorem} {\bf (Haar)} 
Every compact abelian topological group $G$ admits a Haar measure which is invariant under right (and left) translations. This measure is unique up to multiplication by a scalar.
\end{theorem}

\begin{remark} 
$\mu$ is a right Haar measure if and only if
\[ \int_G \Ra f d\mu = \int_G fd\mu, \qquad \text{for all}\; f\in \Co(G), \alpha\in G. \] 
Here, the action of $G$ by translations on the space of complex valued continuous functions $\Co(G)$ is given by:
\[ (\Ra f)(z) = f(z+\alpha). \]
\end{remark}

\begin{example}
\begin{enumerate}[a.]
\item If $G=\UC$, the Haar measure on $\UC$ is the usual Lebesgue measure $d\mu=dx$.
\item If $G=\T^n$, the Haar measure on $\T^n$ is the usual Lebesgue measure 
$d\mu=d\mu_1\times \cdots \times d\mu_n$.
\end{enumerate}
\end{example}

\subsubsection[Statistical behavior of orbits of translations]{Statistical behavior of orbits of translations}
\label{behavoir_orbits}

Let $z\in G$ be any point and $U\subset G$ an open subset. The time that the orbit of $z$ spends in $U$ under iterations of $\Ra$ can be measured as:
\[ \abs{\{0\leq k\leq n-1: \Ra^k(z)\in U\}}. \]
($\abs{\cdot}$ denotes the cardinality of the set.) If the limit
\[\lim_{n\to \infty} \frac{1}{n} \abs{\{0\leq k\leq n-1: \Ra^k(z)\in U\}} \]
exists, then it measures the ``average time'' that the orbit of $z$ spends in $U$. 

Observe that
\[ \Ra^k(z)\in U \qquad \Longleftrightarrow \qquad \chi_U(\Ra^k(z))=1, \]
where $\chi_U$ denotes the characteristic function of $U$. The limit
\[ \lim_{n\to \infty} \frac{1}{n} \sum_{k=0}^{n-1} \chi_U(\Ra^k(z)) \] 
is called the \emph{temporal averages} (Birkhoff's averages, or ergodic averages) of the function $\chi_U$. Now the problem of iteration can be stated as:\\

\noindent \textsc{Problem:} To study the convergence of ergodic averages.\\

If $\varphi\in \Co(G)$, the temporal averages is defined  as: 
\[ \lim_{n\to \infty} \frac{1}{n} \sum_{k=0}^{n-1} \varphi(\Ra^k(z)) \] 
and the ptoblem is to study their convergence. This is the content of the celebrated:

\begin{theorem} {\bf (Birkhoff Ergodic Theorem)} 
If $\varphi\in L^1(G)$, then the limit
\[ \lim_{n\to \infty} \frac{1}{n} \sum_{k=0}^{n-1} \varphi(\Ra^k(z)) \] 
exists for almost all points $z\in G$ (with respect to the Haar measure $\mu$).
\end{theorem}

There are many extensions and several ways to prove the ergodic theorem. Here is a version which uses the notion of invariance.

\begin{definition} 
Let $\varphi\in \Co(G)$. The function $\varphi$ is \emph{invariant} under $\Ra$ if
\[ \varphi(\Ra(z)) = \varphi(z) = \varphi(\Ra^{-1}(z))\qquad \text{for all}\; z\in G. \]
\end{definition}

According with this definition, a function $\varphi$ is invariant under $\Ra$ if $\varphi$ is constant along the orbits of points in $G$ under iterations of $\Ra$.

A subset $U\subset G$ is called \emph{invariant} under $\Ra$ if $\chi_U$ is invariant under 
$\Ra$. That is,
\[ \Ra(U) = U = \Ra^{-1}(U). \]

The translation $\Ra$ is said to be \emph{ergodic} with respect to the Haar measure $\mu$, if for any invariant subset $U$ of $G$ it holds
\[ \mu(U)=0\qquad \text{or}\qquad \mu(U)=1. \]

According with the ergodic theorem, if $\Ra$ is ergodic with respect to the Haar measure 
$\mu$, then
$$ 
\lim_{n\to \infty} \frac{1}{n} \sum_{k=0}^{n-1} \varphi(\Ra^k(z)) = 
\mu(\varphi) := \int_G \varphi d\mu
$$ 
is constant for almost every point $z\in G$ (with respect to $\mu$).

If $\varphi=\chi_U$, then
\[ \lim_{n\to \infty} \frac{1}{n} \sum_{k=0}^{n-1} \chi_U(\Ra^k(z)) = \mu(U). \] 

That is, the orbit of any point in $G$ go through any subset $U$ of positive measure. Moreover, the orbit of any point $z\in G$, under iterations of $\Ra$, remains in the set 
$U$ during an interval of time proportional to $\mu(U)$.

The following remark is an important characterization of ergodicity:

\begin{remark} 
$\Ra$ is ergodic with respect to $\mu$ if and only if the only invariant functions under  
$\Ra$ are constant functions. 
\end{remark}

\subsubsection[Ergodicity of translations on compact abelian groups]{Ergodicity of translations on compact abelian groups}
\label{ergodicity_translations}

The properties shown in previous sections are also valid for translations on any compact abelian group. Recall the following notions from classical Fourier analysis. 

For any function $f\in L^2(G)$, the Fourier coefficients of $f$ are given by
\[ \Ff(\chi):= \int_G f\cdot \overline{\chi} d\mu \qquad (\chi\in \Char(G)). \] 

Since $G$ is compact, it follows that $\Char(G)$ is discrete (countable) and the Fourier series of $f$ (convergent with respect to the norm in $L^2(G)$) can be written as:
\[ f(z) = \sum_{\chi\in \Char(G)} \Ff(\chi)\cdot \chi(z). \]

\begin{lemma}
If $G$ is a monothetic group with generator $\alpha$, then $\chi(\alpha)\neq 1$, for any nontrivial character $\chi\in \Char(G)$.
\end{lemma}

\bproof 
Suppose there exists $\chi\in \Char(G)\setminus \{1\}$ such that $\chi(\alpha) = 1$. Since 
$\chi$ is a homomorpism, it follows that
\[ \chi(n\alpha) = \chi(\alpha)^n = 1, \] 
for any $n\in \Z$. This implies the translates of $\alpha$, $T(\alpha)$ is a proper subset of the kernel of $\chi$, i.e. $T(\alpha)\subset \ker(\chi)$. Since $\chi \neq 1$, it follows that 
$\ker(\chi)$ is a closed and proper subgroup of $G$. This implies that $T(\alpha)$ is not  dense in $G$, that is to say, $G$ is not monothetic.
\eproof

\begin{theorem} 
If $\alpha\in G$, then the following conditions are equivalent:
\begin{enumerate}[a.]
\item The translation $\Ra:G\To G$ is ergodic with respect to the Haar measure.
\item $G$ is a monothetic group with generator $\alpha$.
\end{enumerate}
\end{theorem}

\bproof 
The proof of (a) $\Longrightarrow$ (b) is completely analogous to the cases of the circle and the torus. We prove the implication (b) $\Longrightarrow$ (a).

According to the last Lemma, if $G$ is a monothetic group with generator $\alpha$, 
$\chi(\alpha)\neq 1$, for any nontrivial character $\chi\in \Char(G)$. We prove that $\Ra$ is ergodic with respect to the Haar measure.

Let $f\in L^2(G)$ be an invariant function under $\Ra$. Writing its Fourier series as 
\[ f(z) = \sum_{\chi\in \Char(G)} \Ff(\chi)\cdot \chi(z), \]
this implies
\begin{align*}
f(\Ra(z))	&= \sum_{\chi\in \Char(G)} \Ff(\chi)\cdot \chi(z+\alpha)\\
				&= \sum_{\chi\in \Char(G)} [\Ff(\chi) \chi(\alpha)]\cdot \chi(z).
\end{align*} 

This means that the translation ``acts'' in functions on $L^2(G)$ by multiplication of the Fourier coefficient $\Ff(\chi)$ by $\chi(\alpha)$.

Since $f$ is invariant under $\Ra$, it follows that $f(\Ra(z))=f(z)$. By uniqueness of the Fourier coefficients we conclude that 
\[ \Ff(\chi) = \Ff(\chi) \chi(\alpha). \] 

By hipothesis, $\chi(\alpha)\neq 1$ if $\chi\neq 1$. Hence
\[ \Ff(\chi)=0 \qquad \text{if}\; \chi\neq 1. \]

This implies that $f(z)=\Ff(1)=\text{constant}$. Therefore $\Ra$ is ergodic with respect to 
$\mu$.
\eproof

\section[Basic examples: the circle, the torus and the solenoid]{Basic examples: the circle, the torus and the solenoid}
\label{basic_examples}

This section briefly describes the theory on three fundamental examples of compact abelian groups: the unit circle $\UC$, the torus $\T^2$ and the solenoid $\Ss$. The theory on the first two groups has been developed during the last century. The theory on $\Ss$ is developed very recently in the article \cite{CV}. The general theory presented in Section \ref{rotation_theory} is based on this work. 

\subsection[Rotation sets and dynamics on $\UC$]{Rotation sets and dynamics on $\UC$}
\label{dynamics_circle}

Recall that the unit circle $\UC$ is a compact connected abelian topological group. Denote by $\Homeo_+(\UC)$ the group consisting of all homeomorphisms of $\UC$ which are isotopic to the identity. The suspension $\Sigma_f(\UC)$ of any $f\in\Homeo_{+}(\UC)$ is homeomorphic to $\UC\times\UC$ and
\[ \Hom(\Char(\UC),\R)\cong \Hom(\Z,\R)\cong \R. \]

So the translation set is a subset of $\R$ and consequently the rotation set $\rho(f)$ is a subset of $\UC$. This set consists of a single point and using Theorem \ref{homrot} and the fact that both definitions depends on the integral of the displacement function $f-\id$ with respect to an invariant measure, this number coincides with the usual rotation number of Poincar\'e.

The relationship between the rotation number $\rho(f)$ and the dynamics of the homeomorphisms $f$ is subsumed in the following: 

\begin{theorem} {\bf (Poincar\'e)}
\label{poincare}
Suppose $f\in \Homeo_+(\UC)$ with rotation number $\rho(f)$.
\begin{enumerate}[a.]
\item $\rho(f)$ is rational if and only if $f$ has a periodic point.

\item If $\rho(f)$ is irrational then $f$ is semiconjugated to the rotation $R_{\rho(f)}$. 
The semiconjugation is a conjugation if and only if all the orbits of $f$ are dense.
\end{enumerate}
\end{theorem}

\subsection[Rotation sets and dynamics on $\T^2$]{Rotation sets and dynamics on $\T^2$}
\label{dynamics_torus}

Denote by $\Homeo_+(\T^2)$ the group consisting of all homeomorphisms of $\T^2$ which are isotopic to the identity. The suspension $\Sigma_f(\T^2)$ of any $f\in\Homeo_+(\T^2)$ is homeomorphic to $\T^2\times \UC$ and
\[ \Hom(\Char(\T^2),\R)\cong \R\times \R, \]
since $\Char(\T^2)\cong \Z\times\Z$. 
 
So the translation set is a subset of $\R\times \R$ and consequently the rotation set $\rho(f)$ is a subset of $\T^2$. In the article \cite{MZ}, the authors give some properties of this set. One of the most important is that this set is compact and convex for any 
$f\in \Homeo_+(\T^2)$. The next theorem (see \cite{Fra}) establishes conditions on the rotation set to get periodic orbits as in Theorem \ref{poincare}.
 
\begin{theorem} {\bf(Franks)}
\label{franks}
Let $f\in\Homeo_+(\T^2)$ such that $\rho(f)$ has nonempty interior. If a vector $v$ lies in the interior of $\rho(f)$ and both coordinates are rational, then there is a periodic point 
$x\in \T^n$ with the property that
\[ \frac{F^{q}(x_0) - x_0}{q} = v, \]
where $x_0\in \R^2$ is any lift of $x$, $F$ is a lift of $f$ and $q$ is the least period of 
$x$.
\end{theorem}
 
Later on, in \cite{Jag} the author obtained a version of the semiconjugation theorem as in Theorem \ref{poincare} for irrational pseudo-rotations with bounded mean variation. Here again, the bounded mean variation condition plays a fundamental role.
 
\begin{definition}
$f\in \Homeo_+(\T^2)$ is an \emph{irrational pseudo-rotation} if there exist a vector 
$\rho\in \R^2$ with both coordinates irrational and a lift $F:\R^2\To\R^2$ such that:
\[ \lim_{n\to \infty} \frac{D_{x}(F^n)}{n} = \rho, \qquad (x\in \R^2). \]
\end{definition}

The displacement function $D:\Z\times \R^2\To \R^2$ of a lift $F^n$ of $f^n$ is:
\[ D(n,x) = F^n(x) - x - n\rho. \]

$F$ is said to have \emph{bounded mean variation} with constant $C\in \R$ if 
\[ \norm{ D(n,x)} \leq C \]
for any $n\in \Z$ and $x\in \R^2$.

\begin{theorem} {\bf (J\"ager)}
\label{jager}
Let $f\in \Homeo_+(\T^2)$ be a pseudo-rotation with rotation vector $\rho(f)\in \R^2$ and which satisfies the bounded mean variation condition. 
\begin{enumerate}
\item $\rho(f)$ is totally irrational if and only if $f$ is semiconjugated to $R_{\rho(f)}$.
\item $\rho(f)$ is totally rational if and only  if $f$ has a periodic point.
\end{enumerate}
\end{theorem}

\subsection[Rotation sets and dynamics on $\Ss$]{Rotation sets and dynamics on $\Ss$}
\label{dynamics_solenoid}

In the article \cite{CV} the authors extend the notion of the rotation number for a homeomorphism isotopic to the identity of the one-dimensional universal solenoid $\Ss$. The universal solenoid is a compact connected abelian group and a sort of `generalized circle' described as the projective limit
\[ \Ss = \varprojlim_{n\in \N} \R/n\Z \]
with additional foliated structure. The authors uses duality theory and asymptotic cycles to describe this rotation element. We briefly describe the construction and the details can be consulted in \cite{CV}.

Let $f:\Ss\To \Ss$ be any homeomorphism isotopic to the identity which can be written as $f=\id + \varphi$, where $\varphi:\Ss\To\Ss$ is the displacement function along the one-dimensional leaves of $\Ss$ with respect to the affine structure.

Since $f$ is isotopic to the identity, the suspension space $\SuspS$ of $f$ is homeomorphic to the direct product $\Ss\times \UC$ and again we use the structure of compact abelian topological group of this product. The character group of the suspension is
\[ \Char(\SuspS)\cong \Char(\Ss)\times \Char(\UC)\cong \Q\times \Z. \]

The associated suspension flow $\phi_t:\SuspS\To\SuspS$ is given by:
\[ \phi_t(z,s) = (f^m(z),t+s-m),\qquad (m\leq t+s < m+1). \]

Now, for any given character $\chi_{q,n}\in \Char(\SuspS)$, there exists a unique $1$\nobreakdash-cocycle 
\[ C_{\chi_{q,n}}:\R\times \SuspS\To \R \] 
associated to $\chi_{q,n}$ such that
\[ \chi_{q,n}(\phi_t(z,s)) = \exp(2\pi iC_{\chi_{q,n}}(t,(z,s)))\cdot \chi_{q,n}(z,s), \] 
for every $(z,s)\in \SuspS$ and $t\in \R$. From here it is obtained an explicit expression for the $1$\nobreakdash-cocycle $C_{\chi_{q,n}}(t,(z,s))$ and, by Birkhoff's Ergodic Theorem, there is a well-defined homomorphism
\[ H_f:\Char(\SuspS)\To \R \] 
given by
\[ H_f(\chi_{q,n}) = \int_{\SuspS} C_{\chi_{q,n}}(1,(z,s)) d\nu, \] 
where $\nu$ is a $\phi_t$-invariant Borel probability measure on $\SuspS$. The explicit calculation of the $1$\nobreakdash-cocycle implies that the last integral only depends on an $f$-invariant Borel probability measure $\mu$. So, if $\PP_f(\Ss)$ is the weak$^*$ compact convex space consisting of all such measure on $\Ss$, the well-defined continuous homomorphism
\[ \rho_\mu(f):\Char(\Susp)\To \UC \] 
given by
\[ \rho_\mu(f)(\chi_{q,n}) := \exp(2\pi iH_f(\chi_{q,n})) \] 
determines an element in $\Char(\Char(\Susp))\cong \Ss\times \UC$ not depending on the second component. By Pontryagin's duality theorem, it defines an element 
$\rho_\mu(f)\in \Ss$ called the \emph{rotation element} associated to $f$, which is the generalized Poincar\'e rotation number.

As expected, $\rho_\mu(f)$ is an element in the solenoid itself and it measures, in some sense, the average displacement of points under iteration of $f$ along the one-dimensional leaves with the Euclidean metric.

If $\rho:\PP_f(\Ss)\To \Ss$ is the map given by $\mu\mTo \rho_\mu(f)$, then $\rho$ is continuous from $\PP_f(\Ss)$ to $\Ss$. Since $f$ is isotopic to the identity, the image 
$\rho(\PP_f(\Ss))$ is a compact subset of $\Ss$.

\begin{definition} 
The \emph{rotation set} of $f$ is $\rho(\PP_f(\Ss))$.
\end{definition}

Since the pathwise components of $\Ss$ are one-dimensional leaves, this set 
$\rho(\PP_f(\Ss))$ is a compact interval $I_f$, which, up to a translation, is contained in the one-parameter subgroup $\Ll_0$. Since $\Ll_0$ is canonically isomorphic to $\R$, it is possible to identify $I_f$ with an interval in the real line. In particular, if $f$ is uniquely ergodic, then the interval $I_f$ reduces to a point and in this case, this unique element of 
$\Ss$ is the rotation element of $f$.

Let us recall the following definition (\cite{CV}):

\begin{definition}
A homeomorphism $f:\Ss\To \Ss$ is said to have \emph{bounded mean variation} if there exists $C>0$ such that the sequence 
\[ \{F^n(z) - z - n\tau(F)\}_{n\geq 1} \] 
is uniformly bounded by $C$. Here, $F$ is any lift of $f$ to $\R\times \Zz$, $\tau(F)$ is a lifting of $\rho(f)$ to $\R\times \Zz$ and $z\in \R\times \Zz$.
\end{definition}

This definition does not depend on the lifting. The first part of the generalized Poincar\'e Theorem in \cite{CV} is stated as follows:

\begin{theorem}{\bf(Cruz-L\'opez and Verjovsky)}
\label{cv} 
Suppose that $f:\Ss\To \Ss$ is any homeomorphism isotopic to the identity with irrational rotation element $\rho(f)\in \Ss$. The homeomorphism $f$ is semiconjugated to the irrational rotation $R_{\rho(f)}$ if and only if $f$ has bounded mean variation.
\end{theorem}

The question remains:

\begin{question}
\label{question_conjugacy}
Under the same hypothesis of this theorem, is $f$ conjugated to the rotation $R_{\rho(f)}$ when $f$ is minimal?
\end{question}

\section[Entropy and rotation sets]{Entropy and rotation sets}
\label{entropy_rotation}

This section presents some relations between rotation sets and topological entropy (that will be called just entropy) of homeomorphisms of certain compact abelian groups. Denote the entropy of a homeomorphism $f$ of a topological space by $h_{top}(f)$. Entropy is a useful tool to estimate the complexity of a given dynamical system, which can be seen as the growth rate of the number of orbits of the system. A general reference on this topic is \cite{KH} where it can be found the explicit definition of entropy. 

\begin{example}
For any $f\in\Homeo_+(\UC)$ the entropy $h_{top}(f)$ is 0 (see \cite{KH}, Corollary 11.2.10). It should be emphasized that the behavior of orbits of points under iteration of $f$ is determined by the rotation number $\rho(f)$. Actually, this is true for any homeomorphism of the circle.
\end{example}

In the case of the universal solenoid a similar result holds.

\begin{example}
For any $f\in\Homeo_+(\Ss)$, $h_{top}(f)=0$ (see \cite{Kwa2}, Theorem 2). In this case we do not know if the behavior of orbits is only determined by the rotation set, but there are some works in this line cited before.

However, it is not true that any homeomorphism of the solenoid has null  entropy. In \cite{Kwa2} it is proved that the topological entropy of these homeomorphisms only depends on the isotopy class. For instance, multiplication by 2, $q\overset{f_2}{\mTo 2q}$, is an automorphism of $\Q$ and, since $\Ss$ is the Pontryagin dual of $\Q$, it induces a dual map $\widehat{f}_2:\Ss\To \Ss$. This map is an expanding map and has positive entropy.
\end{example}

The case of the torus $\T^2$ is more subtle because it is not true that for any 
$f\in\Homeo_+(\T^2)$, $h_{top}(f)=0$. In \cite{LlM} it is proved the following theorem:

\begin{theorem}
If $f\in\Homeo_+(\T^2)$ has three or more periodic orbits with noncollinear rotation vectors, then $f$ has positive entropy.
\end{theorem}

In particular this implies, using Theorem \ref{franks}, that for any $f\in\Homeo_+(\T^2)$ such that $\rho(f)$ has nonempty interior, the entropy is positive. The authors also prove that $f$ has a kind of chaos called ``toroidal chaos''.

In this line of ideas, in \cite{Kwa1} appears the definition of \emph{topological entropy at the rotation vector} for any $f\in\Homeo_+(\T^2)$. Topological entropy at the rotation vector is a function with domain of definition the rotation set $\rho(f)$ taking real values. The value of the topological entropy of $f$ is its maximum value. The author shows that the topological entropy is positive at the rotation vector that lies in the interior of the rotation set. This idea has been generalized very recently in \cite{KW1} and \cite{KW2}.

\begin{remark}
As a conclusion of this section the following questions arise: 
\begin{enumerate}
\item What is the relationship between entropy and rotation sets for general compact abelian groups?
\item What is the relation between the rotation set and the set of periodic orbits?
\item Can be generalized the definition of topological entropy at the rotation vector for getting similar results as in \cite{Kwa1}?
\end{enumerate}
\end{remark}

Strongly related to Question (2) is the analysis of rationality for rotation sets of homeomorphisms of solenoids made in \cite{Lop}.

\section[Final remarks]{Final remarks}
\label{final_remarks}

In order to be able to describe dynamical properties of a general homeomorphim of a compact abelian group, it is important to understand the structure of the complete group of homeomorphisms $\Homeo(G)$ of $G$. The first kind of results is already provided in \cite{Kee}:
\[ \Homeo(G)\simeq G\times \Homeo_e(G), \]
where $G$ is identified with the group of translations and $\Homeo_e(G)$ is the subgroup containing all homeomorphisms which leaves fixed the identity element $e\in G$. When 
$G=\Ss$ is a one-dimensional solenoidal group, the author proved much more: 
\[ \Homeo_e(\Ss)\simeq \Aut(\Ss)\times \Homeo_0(\Ss) \]
where $\Aut(\Ss)$ denotes the group of continuous automorphisms of $\Ss$ and 
$\Homeo_0(\Ss)$ denotes the set of all homeomorphisms isotopic to the identity. In this case the group $\Homeo_0(\Ss)$ is contractible.

So the first question arises: 
\begin{question}
For a general compact connected abelian group $G$, is it possible to get a decomposition 
$\Homeo_e(G)\simeq \Aut(G)\times \Homeo_0(G)$ with $\Homeo_0(G)$ contractible?
\end{question} 
 
In the case of an affirmative answer for this question, one could analyze the dynamics of the particular cases $f\in \Aut(G)$ and $f\in\Homeo_0(G)$. The case of translations (which are contained in $\Homeo_0(G$)) was already studied in Section \ref{dynamics_translations} where it was established that any translation $\Ra$ is ergodic if and only if the group $G$ is monothetic with monothetic generator $\alpha$. For $\Aut(G)$, the next theorem in \cite{Hal} does the job:

\begin{theorem} {\bf (Halmos)}
A continuous automorphism of a compact abelian group $G$ is ergodic if and only if the induced automorphism on the character group has no finite orbits.
\end{theorem}

Corresponding to the factor $\Homeo_0(G)$ and according with the ideas mentioned in this note, for any compact abelian topological group $G$ there is a well-defined rotation set. However, it can be trivial, or perhaps it couldn't give us any dynamical information. Even in the case when the theory is not trivial we do not know general properties of the rotation set, for example, it is unknown to us if the rotation set is always compact or closed. Also, one can ask for the convexity of the rotation set when $\Hom(\Char(G),\R)$ is equipped with the structure of an $\R$-vector space.

The natural question, according with classical Rotation Theory, states:

\begin{question}
Does there exist a semiconjugation theorem with the correct hypothesis as in Section \ref{dynamics_solenoid}?
\end{question}

\section*{Acknowledgements}

The authors gratefully acknowledge Cynthia Verjovsky Marcotte for extensive editing of this paper. 
The first two authors are deeply grateful with the Department of Mathematics of the University of El Salvador and the Ministry of Education of El Salvador for providing an excellent atmosphere which was essential in developing this paper. Special thanks to Ing. Carlos Canjura, former Minister of Education and Dr. Nerys Funes, former Head of the Department for the financial support and their encouragement while working on the project.


\baselineskip=17pt

\begin{thebibliography}{99}

\bibitem{AJ} Aliste-Prieto, J. and J\"ager, T. \emph{Almost periodic structures and the semiconjugacy problem}. Jour. Diff. Eq. 252, (2012), 4988--5001.  

\bibitem{Ath} Athanassopoulos, K. \emph{Some aspects of the theory of asymptotic cycles}. Expo. Math. 13 (1995), 321--336.

\bibitem{CC} Candel, A. and Conlon, L. \emph{Foliations I}. American Mathematical Society, 2000.

\bibitem{CFS} Cornfeld, I.P., Fomin, S.V. and Sinai, Y.G. \emph{Ergodic Theory}, Springer-Verlag, 1982.

\bibitem{Cro} Crovisier, S. \emph{Exotic rotations}. In \'Ecole d'\'et\'e M\'ethodes topologiques en dynamique des surfaces, Grenoble, (2006).

\bibitem{CV} Cruz-L\'opez, M. and  Verjorvsky, A. \emph{Poincar\'e theory for the 
ad\`ele class group $\mathbb{A}/\Q$ and compact abelian one-dimensional solenoidal groups.} Submitted for publishing.

\bibitem{Fra} Franks, J. \emph{Realizing rotation vector for torus homeomorphisms}. Transactions of the American Mathematical Society 311 (1989), 107--115.

\bibitem{FHK} Furno, J., Haynes, A. and Koivusalo, H. \emph{Bounded remainder sets for rotations on the adelic torus}. Proc. Amer. Math. Soc. 147 (2019), no. 12, 5105--5115. 

\bibitem{Ghys} Ghys E. \emph{Groups acting on the circle}. L$'$Enseignement Math\'ematique 47 (2001), 329--407.

\bibitem{Hal} Halmos, P.R. \emph{On automorphisms of compact groups}.  Bull. Amer. Math. Soc. 49 (1943), 619--624.

\bibitem{Her} Herman, M.R. \emph{Sur la conjugaison diff\'erentiable des diff\'eomorphismes du cercle \`a des rotations}. Publ. Math. IHES 49 (1979), 5--233.

\bibitem{HR} Hewitt, E. and Ross, K.A. \emph{Abstract Harmonic Analysis}, Springer-Verlag, 1979.

\bibitem{Jag} J\"ager, T. \emph{Linearization of conservative toral homeomorphisms}. Invent. Math. 3 (2009),  pp. 601--616.

\bibitem{Kat} Katok, A. \emph{Fifty years of entropy in dynamics: 1958–2007.} J. Mod. Dyn. 1 (2007), 545--596.

\bibitem{KH} Katok, A. and Hasselblatt, B. \emph{Introduction to the Modern Theory of Dynamical Systems}, Cambridge University Press, 1997.

\bibitem{Kee} Keesling, J. \emph{The group of homeomorphisms of a solenoid}. Trans. Amer.Math. Soc. 172 (1972), 119--131.

\bibitem{KW1} Kucherenko, T. and Wolf, Ch. \emph{Entropy and rotation sets: a toy model approach}. Commun. Contemp. Math. 18 (2016), 1550083, 23 pp. 

\bibitem{KW2} Kucherenko, T. and Wolf, Ch. \emph{Geometry and entropy of generalized rotation sets}. Israel J. Math. 199 (2014), 791--829.

\bibitem{Kwa1} Kwapisz, J.M. \emph{Rotation Sets and Entropy}. PhD Dissertation, SUNY at Stony Brook, 1995.

\bibitem{Kwa2}  Kwapisz, J.M. \emph{Homotopy and dynamics for homeomorphisms of solenoids and Knaster continua.} Fundamenta Matematicae 168 (2001), 251--278.

\bibitem{LlM} Llibre, J. and MacKay, R.S. \emph{Rotation vectors and entropy for homeomorphisms of the torus isotopic to the identity}. Ergod. Th. $\&$ Dynam. Sys. 11 (1991) 115--128.

\bibitem{Lop} L\'opez-Hern\'andez, F.J. \emph{Dynamics of induced homeomorphisms of one-dimensional solenoids}. Discrete Contin. Dyn. Syst. 38 (2018), 4243--4257.

\bibitem{Mis} Misiurewicz, M. \emph{Rotation Theory}. Online Proceedings of the RIMS Workshop ``Dynamical Systems and Applications: Recent Progress'', 2006. 

\bibitem{MZ} Misiurewicz, M. and Ziemian, K. \emph{Rotation sets for maps of tori}. J. London Math. Soc.  3 (1998) pp. 490--506.

\bibitem{Poi} Poincar\'e, J.H. \emph{Memoire sur les courbes d\'efinis par une \'equation diff\'erentielle}. Journal de Math\'ematiques. 7 (1881), 375--422.

\bibitem{Pol} Pollicott, M. \emph{Rotation sets for homeomorphisms and homology}. Trans. Amer. Math. Soc. 331  (1992), 881--894.

\bibitem{Pon} Pontryagin, L. \emph{Topological groups}. Gordon and Breach, 1966.

\bibitem{Sch} Schwartzman, S. \emph{Asymptotic cycles}, Ann. of Math. 66  (1957), 270--284.

\bibitem{Ste} Steenrod, N. \emph{Universal homology groups}, Amer. J. Math. 58 (1936), 661--701.

\bibitem{Wil} Wilson,  J.S. \emph{Profinite groups}, Oxford University Press, 1990.

\end{thebibliography}
\end{document}